\documentstyle[12pt,amscd]{amsart}
\newtheorem{pr}{Proposition}
\newtheorem{lm}{Lemma}

\newcommand{\proj}{\bold P}

\newcommand{\barr}{\overline}

\newcommand{\rarr}{\rightarrow}
\newcommand{\oh}{{\cal{O}}}
\newcommand{\M}{\barr{M}_{0,n}(\proj^r,d)}

\newcommand{\eqq}{\stackrel{\sim}{=}}
\newcommand{\deli}{\bigtriangleup}
\newcommand{\com}{\Bbb{C}}
\newcommand{\LL}{\cal{L}}
\newcommand{\HH}{\cal{H}}
\newcommand{\DD}{\cal{D}}
\newcommand{\KK}{\cal{K}_{\barr{M}}}
\newcommand{\NN}{\cal{N}}
\begin{document}
\title{The Canonical Class of $\barr{M}_{0,n}(\proj^r, d)$
and Enumerative Geometry}
\author{Rahul Pandharipande}
\date{5 September 1995}
\maketitle
\pagestyle{plain}
\setcounter{section}{-1}
\section{Summary}
Let $\com$ be the field of complex numbers.
Let
the Severi variety $$S(0,d)\subset
\proj\big(H^0(\proj^2, \oh_{\proj^2}(d))\big)$$
be the quasi-projective locus of irreducible, nodal
rational curves. Let $\barr{S}(0,d)$ denote the closure
of $S(0,d)$.
Let $p_1, \ldots, p_{3d-2}$ be $3d-2$ general points in $\proj^2$.
Consider the subvariety $C_d\subset \barr{S}(0,d)$ corresponding
to curves passing through all the  points $p_1, \ldots, p_{3d-2}$.
$C_d$ is a complete curve in $\proj
\big(H^0(\proj^2, \oh_{\proj^2}(d))\big)$.
Let $N_d$ be the degree of $C_d$. $N_d$ is
determined by the recursive relation
([K-M], [R-T]):
$$N_1=1$$
$$\forall d>1, \ \ \ N_d= \sum_{i+j=d, \ i,j>0}
N_iN_j \bigg( i^2j^2 {3d-4 \choose 3i-2} -
i^3j {3d-4 \choose 3i-1} \bigg).$$
For $d\geq 3$, $C_d$ is singular. The
arithmetic genus $g_d$ of $C_d$ is determined by:
$$g_1=0,$$
$$g_2=0,$$
$$2g_d-2=
{6d^2+5d-15\over 2d}N_d +
{1\over 4d}\sum_{i=1}^{d-1} N_i N_{d-i}\Big(15i^2(d-i)^2
-8di(d-i)-4d\Big) {3d-2 \choose 3i-1}.$$
The last formula holds for $d\geq 3$.
The geometric genus $\tilde{g}_{d}$ of the
normalization $\tilde{C}_d$ is determined by ($d\geq 1$):
$$2\tilde{g}_{d}-2 = -{3d^2-3d+4\over 2d^2} N_d
+ {1\over 4d^2} \sum_{i=1}^{d-1} N_i N_{d-i}
(id-i^2) \Big( (9d+4)i(d-i)-6d^2\Big) {3d-2\choose 3i-1}.$$
These genus formulas
are established by adjunction and intersection
on Kontsevich's space of stable maps $\M$.

The author thanks D. Abramovich for many useful
remarks and for finding an error
in the singularity analysis in an earlier version of
this paper.

\section {The Canonical Class of $\M$}
Let $\M$ be the coarse moduli space of degree $d$, Kontsevich
stable maps from $n$-pointed, genus $0$ curves to
$\proj^r$. Foundational treatments of $\M$ can be found
in [Al], [P1], [K], and [K-M]. Only the case $r\geq 2$
will be considered here.
Let $\LL_p$ denote the line bundle obtained on $\M$
by the $p^{th}$ evaluation map $(1\leq p \leq n)$.
Let $\deli$ be the set of boundary divisors.
Let $\HH$ denote the divisor of maps meeting
a fixed codimension 2 linear space of $\proj^r$.
$\HH=0$ if $d=0$.
In [P2], it is shown the
classes
$$\{\LL_p\} \cup \deli \cup \{\HH\}$$
span $Pic(\M) \otimes \Bbb{Q}$.

The canonical class of the stack $\barr{\cal{M}}_{0,n}(\proj^r,d)$
has the following coarse moduli interpretation.
$\M$ is an irreducible variety with finite quotient
singularities. When $r\geq 2$,
the automorphism-free locus $\barr{M}^*_{0,n}(\proj^r,d)
\subset \M$ is a nonsingular,
fine moduli space with codimension 2 complement except
when ([P2]) $$[0,n,r,d]=[0,0,2,2].$$ For $r\geq 2$ and
$[0,n,r,d]\neq [0,0,2,2]$, the first Chern class of the
cotangent bundle to the moduli space $\barr{M}^*_{0,n}(\proj^r,d)$
yields the canonical class in $Pic(\M) \otimes \Bbb{Q}$.

Let $P=\{1,2,\ldots, n\}$ be the set of markings ($P$ may be the
empty set).
The boundary components
are in bijective correspondence with
data of weighted partitions $(A\cup B, d_A, d_B)$ where:
\begin{enumerate}
\item[(i.)] $A\cup B$ is a partition of $P$.
\item[(ii.)] $d_A+d_B=d$, $d_A \geq 0$, $d_B \geq 0$.
\item[(iii.)] If $d_A=0$ (resp. $d_B=0$), then $|A|\geq 2$ (resp.
$|B| \geq 2$).
\end{enumerate}
Define $\DD_{i,j}$ to be the reduced sum of boundary
components with $d_A=i$ and $|A|=j$.
Note $0\leq i \leq d$ and $0\leq j \leq n$.
The divisors $\DD_{0,0}$, $\DD_{0,1}$, $\DD_{d,n-1}$, $\DD_{d,n}$
are equal to 0 by stability. Also, $\DD_{i,j}=\DD_{d-i,n-j}$.

Consider first the case $d=0$. $\barr{M}_{0,n}(\proj^r,0) \eqq
\barr{M}_{0,n} \times \proj^r$.
It suffices to determine the canonical class of $\barr{M}_{0,n}$.
\begin{pr} The canonical class $\KK$ of $\barr{M}_{0,n}$ is
determined in $Pic(\barr{M}_{0,n}) \otimes \Bbb{Q}$ by:
\begin{equation}
\label{mumkn}
\KK =
\sum_{j=2}^{[{n\over 2}]} \Big( {j(n-j) \over n-1} -2 \Big) \DD_{0,j}.
\end{equation}
\end{pr}
\noindent
The canonical class has a different form in case
$d>0$, $n=0$, $r\geq 2$.
\begin{pr} The canonical class $\KK$ of
$\barr{M}_{0,0}(\proj^r, d)$ ($d>0, r\geq 2$) is determined in
$Pic(\barr{M}_{0,0}(\proj^r,d) \otimes
\Bbb{Q}$ by:
\begin{equation}
\label{genny}
\KK= -{(d+1)(r+1)\over 2d} \HH + \sum_{i=1}^{[{d\over 2}]}
\Big( {(r+1)(d-i)i\over 2d}-2 \Big) \DD_{i,0}.
\end{equation}
\end{pr}
\noindent Finally, when $d>0$, $n>0$, $r\geq 2$, the form of the
canonical class is the following:
\begin{pr} The canonical class of $\KK$ of
$\M$ ($d>0$, $n>0$, $r\geq 2$) is determined in
$Pic(\M) \otimes \Bbb{Q}$ by:
\begin{equation}
\label{fini}
\KK=  -  {(d+1)(r+1)d-2n \over 2d^2} \HH
 - \sum_{p=1}^{n} {2\over d}\LL_p
\end{equation}
$$+\sum_{i=0}^{[{d\over 2}]} \sum_{j=0}^{n}
\Big({(r+1)(d-i)di + 2d^2j-4dij+2ni^2 \over 2d^2}-2 \Big) \DD_{i,j}.$$
\end{pr}
\noindent
Equation (\ref{mumkn}) can be derived from the explicit
construction of $\barr{M}_{0,n}$ described in [F-M].
Equations (\ref{mumkn}-\ref{fini})
will be established here via intersections with curves.

\section{Computing The Canonical Class}
\subsection{Curves in $\M$}
By Proposition (2) of [P2], the canonical projection
$$Pic(\M)\otimes \Bbb{Q} \rarr Num(\M) \otimes \Bbb{Q}$$ is
 an isomorphism. Hence, the canonical class of $\M$ can
be established via intersections with curves.
Curves can easily be found in $\M$. The notation of [P2]
is recalled here.

Let $C$ be a nonsingular, projective curve.
Let $\pi: S=\proj^1 \times C \rarr C$. Select $n$ sections
$s_1, \ldots, s_n$
of $\pi$.
A point $x\in S$ is an
{\em intersection point} if two or more sections contain $x$.
Let $\NN$ be a line bundle on $S$ of type $(d,k)$.
Let $z_l\in H^0(S, \NN)$  $(0\leq l \leq r)$
determine a rational map $\mu: S - \rarr \proj^r$ with simple
base points. A point $y\in S$ is a {\em simple base point} of
degree $1\leq e\leq d$ if
the blow-up of $S$ at $y$ resolves $\mu$ locally
at $y$ and the resulting map is of degree $e$ on the exceptional
divisor $E_y$.
The set of {\em special points} of $S$ is the
union of the intersection points and the simple base points.
Three conditions are required:
\begin{enumerate}
\item [(1.)] There is at most one special point in each fiber of
$\pi$.
\item [(2.)] The sections through each intersection point $x$
have distinct tangent directions at $x$.
\item [(3.)]
\begin{enumerate}
\item[(i.)] $d=0$.  No $n-1$ sections pass through a point $x\in S$.
\item[(ii.)] $d>0$.
If at least $n-1$ sections pass through a
point $x\in S$, then $x$ is not a simple base point of degree $d$.
\end{enumerate}
\end{enumerate}
Condition (3.ii) implies  there are no simple base points of degree $d$
if $n=0$ or $1$.
Let $\barr{S}$ be the blow-up of $S$ at the special points.
It is easily seen $\barr{\mu} : \barr{S} \rarr \proj^r$ is
Kontsevich stable family of $n$-pointed, genus $0$ curves over
$C$. Condition (2) ensures the strict transforms of the sections
are disjoint. Condition (3) implies Kontsevich stability.
There is a canonical morphism $\lambda:C \rarr \M$.
Condition (1)
implies $C$ intersects the boundary components transversally.

\subsection{$\barr{M}_{0,n}$}
Curves in $\barr{M}_{0,n}$ are obtained by the above
construction (omitting the map $\mu$).
Let $s_1, \ldots, s_n$ be n sections of $\pi:S=\proj^1 \times C
 \rarr C$ satisfying $(1), (2), (3.i)$.
For $1\leq \alpha \leq n$, let
$s_\alpha$ be of type $(1,\sigma_\alpha)$ on $S=P^1\times C$.
Let $\barr{\pi}:\barr{S}
\rarr C$ be the blow-up of  $S$ as above.
Let $s_1, \ldots, s_n$ also denote the transformed sections of
$\barr{\pi}$. Let $Q$ denote the points of $C$ lying under the
special points of $S$.
There is a canonical sequence
$$0 \rarr R^1\barr{\pi}_*(\omega^*_{\barr{\pi}}(-\sum_{1}^{n} s_\alpha))
\rarr
\lambda^*(T_{\barr{M}_{0,n}})
\rarr \bigoplus_{q\in Q} \Bbb{C}_q \rarr 0.$$
(See, for example, [K].)
Hence $C\cdot \KK=- deg \big(
R^1\barr{\pi}_*(\omega^*_{\barr{\pi}}(-\sum_{1}^{n} s_\alpha)) \big) -
C\cdot \sum_{j=2}^{[{n\over 2}]} \DD_{0,j}$.

The degree of
$R^1\barr{\pi}_*(\omega^*_{\barr{\pi}}(-\sum_{1}^{n} s_\alpha))$
is determined by the Grothendieck-Riemann-Roch formula.
Let $x_j$ for $2 \leq j \leq n-2$ be the number of intersection points
of $S$ which lie on exactly $j$ sections.
If $j\neq n/2$, $C\cdot \DD_{0,j}= x_j+ x_{n-j}$.
If $j=n/2$, $C\cdot \DD_{0,j}
=x_j$.
G-R-R yields:
$$ deg \big(R^1\barr{\pi}_*
(\omega^*_{\barr{\pi}}(-\sum_{1}^{n} s_\alpha)) \big)
=\sum_{1}^{n}2 \sigma_{\alpha} + \sum_{2}^{n-2} (1-j)x_j.$$
By the transverse intersection condition, the following relation
must hold:
$$\sum_{1}^{n} \sigma_{\alpha} =
{1 \over n-1} \sum_{2}^{n-2} {j^2-j\over 2}x_j.$$
Combining equations yields:
$$C\cdot \KK =
\sum_{2}^{n-2} \Big( j-2 - {j^2-j\over n-1} \big) x_j$$
$$ =  \sum_{2}^{n-2} \big({j(n-j)\over n-1} -2 \Big) x_j.$$
Hence both sides of equation (\ref{mumkn}) have the same
intersection numbers with $C$.
Let $D$ be any nonsingular curve in $\barr{M}_{0,n}$ which
intersects the boundary transversely. The universal family over
$D$ can be blown-down to a projective bundle $\pi: T \rarr D$.
The above calculation covers the case where $T=\proj^1\times D$.
The general case (in which T is {\em any} $\proj^1$-bundle) is
identical. Since $A^1(\barr{M}_{0,n})$ is spanned by curves
meeting the boundary transversely,
Proposition (\ref{mumkn}) is immediate.

\subsection{$\barr{M}_{0,0}(\proj^r,d)$}
The case $d>0$, $n=0$, $r\geq 2$ is now considered.
Let $\barr{\pi}:\barr{S} \rarr C$,
$\barr{\mu}:\barr{S} \rarr \proj$
be a family of stable maps as above. There is  canonical
exact sequence
$$0 \rarr \barr{\pi}_*(\omega^*_{\barr{\pi}}) \rarr
 \barr{\pi}_*\barr{\mu}^*(T_{\proj^r})
\rarr \lambda^*(T_{\barr{M}_{0,0}
(\proj^r,d)}) \rarr \bigoplus_{p\in Q} \Bbb{C}_p \rarr 0.$$
Hence $C\cdot \KK=
+deg \big( \barr{\pi}_*(\omega^*_{\barr{\pi}}) \big)
-deg \big(\barr{\pi}_*\barr{\mu}^*(T_{\proj^r}))\big)
          -\sum_{1}^{[{d\over 2}]} \DD_{i,0}$.
Let $x_i$
for $1\leq i \leq d-1$ be the number of simple base points
of $\mu: S - \rarr \proj^r$ of degree exactly $i$.
If $i\neq d/2$, $C\cdot \DD_{i,0}= x_i + x_{d-i}$. If
 $i=d/2$, $C \cdot \DD_{i.0}=x_i$.
Via G-R-R,
$$deg\big( \barr{\pi}_*(\omega^*_{\barr{\pi}}) \big) =
-\sum_{1}^{d-1} x_i,$$
$$deg\big( \barr{\pi}_*(\barr{\mu}^*(T_{\proj^r})) \big)=
(r+1)(d+1)k - \sum_{1}^{d-1} {(r+1)(i^2+i)\over 2} x_i.$$
Combining equations yields:
$$C\cdot \KK =  -(r+1)(d+1)k + \sum_{1}^{d-1}
\Big( {(r+1)(i^2+i)\over 2}-2 \Big) x_i.$$
Finally $C\cdot \HH$ must be computed:
$$C \cdot \HH = 2dk- \sum_{1}^{d-1} i^2 x_i.$$
These equations (plus algebra) verify  Proposition (2).
As before, the complete proof requires the above calculation
for {\em any} $\proj^1$-bundle $\pi: S \rarr C$. Again, the
generalization to this case is trivial.

\subsection{$\M$}
Finally, the case $d>0$, $n>0$, $r\geq 2$ is considered.
Let $\barr{\pi}:\barr{S} \rarr C$, $\barr{\mu}:\barr{S} \rarr \proj$
be a family of stable maps as above.
Let $s_1, \ldots, s_n$ be n sections of $\pi:S=\proj^1 \times C
 \rarr C$ satisfying $(1), (2), (3.ii)$.
For $1\leq \alpha \leq n$, let
$s_\alpha$ be of type $(1,\sigma_\alpha)$ on $S=P^1\times C$.
There is a canonical exact sequence
$$0 \rarr \barr{\pi}_*(\omega^*_{\barr{\pi}}) \rarr
\barr{\pi}_*(\omega^*_{\barr{\pi}}|_{\sum s_{\alpha}}) \bigoplus
 \barr{\pi}_*\barr{\mu}^*(T_{\proj^r})
\rarr \lambda^*(T_{\barr{M}_{0,0}
(\proj^r,d)}) \rarr \bigoplus_{p\in Q} \Bbb{C}_p \rarr 0.$$
Hence $C\cdot \KK=
+deg \big( \barr{\pi}_*(\omega^*_{\barr{\pi}}) \big)
-(\omega^*_{\barr{\pi}} \cdot \sum_{1}^{n} s_{\alpha})
-deg \big(\barr{\pi}_*\barr{\mu}^*(T_{\proj^r}))\big)
          -\sum_{i=1}^{[{d\over 2}]} \sum_{j=0}^{n} \DD_{i,j}$.
Let $z_{i,j}$
for $0\leq i \leq d$ and $0\leq j \leq n$
be the number of special points of $S$
that are simple base points of degree exactly $i$ and
that lie on exactly $j$ sections.
If $i\neq d/2$ or $j\neq n/2$, then
 $C \cdot \DD_{i,j}= z_{i,j}+ z_{d-i, n-j}$.
If $i=d/2$ and $j=n/2$, then $C\cdot \DD_{i,j}=z_{i,j}$.
Via G-R-R,
$$deg\big( \barr{\pi}_*(\omega^*_{\barr{\pi}}) \big) =
-\sum_{i=0}^{d}\sum_{j=0}^{n} z_{i,j},$$
$$deg\big( \barr{\pi}_*(\barr{\mu}^*(T_{\proj^r})) \big)=
(r+1)(d+1)k -
\sum_{i=0}^{d}\sum_{j=0}^{n} {(r+1)(i^2+i)\over 2} z_{i,j}.$$
A simple calculation yields:
$$\omega^*_{\barr{\pi}}
\cdot \sum_{1}^{n} s_{\alpha}= \sum_{1}^{n} 2\sigma_\alpha
- \sum_{i=0}^{d} \sum_{j=0}^{n} j z_{i,j}.$$
The two additional intersection numbers are:
$$C \cdot \HH= 2dk - \sum_{i=0}^{n} \sum_{j=0}^{n} i^2 z_{i,j},$$
$$C \cdot \sum_{1}^{n} \LL_p = nk+ \sum_{1}^{n}d \sigma_{\alpha}
-  \sum_{i=0}^{n} \sum_{j=0}^{n} ij z_{i,j}.$$
Now algebra yields the desired equality of intersections that
establishes Proposition (3). Again the calculation must
be done in case  $\pi: S \rarr C$ is a $\proj^1$ bundle.

\section{The genus of $C_d$, $\tilde{C}_d$}

\subsection{Singularities}
Let $C_d\subset \barr{S}(0,d)$ be the dimension $1$ subvariety
corresponding to curves passing through $3d-2$ general points
$p_1,\ldots, p_{3d-2}$ in $\proj^2$. Let $\hat{C}_d\subset
\barr{M}_{0,0}(\proj^2, d)$ be the dimension $1$ subvariety
corresponding to maps passing through $p_1, \ldots, p_{3d-2}$.
The singularities of $C_d$, $\hat{C}_d$ will be analyzed.

Let $[\mu]\in \barr{M}_{0,0}(\proj^2,d)$ correspond to
an automorphism-free map with domain $\proj^1$. There is
a normal sequence on $\proj^1$:
$$0 \rarr T_{\proj^1}
\stackrel {d\mu}{\rarr}
 \mu^*(T_{\proj^2}) \rarr Norm \rarr 0.$$
The tangent space
to $\barr{M}_{0,0}(\proj^2,d)$ is $H^0(\proj^1, Norm)$.
If $\mu$ is an immersion, $Norm \eqq \oh_{\proj^1}(3d-2)$.
If $\mu$ is not an immersion $Norm$ will have torsion.
A {\em $1$-cuspidal} rational plane curve is
an irreducible rational plane curve with nodal singularities
except for exactly 1 cusp.
If $\mu$ corresponds to a $1$-cuspidal rational curve, then
there is a sequence:
$$0 \rarr \com_p \rarr Norm \rarr \oh_{\proj^1}(3d-3) \rarr 0$$
where $p$ is the point of $\proj^1$ lying over the cusp.
Since $3d-2$ distinct points of $\proj^1$ always
impose independent conditions on
$H^0(\proj^1, \oh_{\proj^1}(3d-2))$ and
$H^0(\proj^1, \oh_{\proj^1}(3d-3))$, Lemma (\ref{t}) has
been established:
\begin{lm}
\label{t}
Let $[\mu]\in \hat{C}_d$ be a point corresponding to an
irreducible, nodal or $1$-cuspidal rational curve
with all singularities
distinct
from $p_1, \ldots, p_{3d-2}$. $\hat{C}_d$ is nonsingular at $[\mu]$.
\end{lm}
\noindent The corresponding analysis for $C_d$ is more involved.
\begin{lm}
\label{tt}
Let $x\in C_d$ be a point corresponding to an
irreducible, nodal rational curve with nodes
distinct
from $p_1, \ldots, p_{3d-2}$. $C_d$ is nonsingular at $x$.
\end{lm}
\begin{pf}
Let $X\subset \proj^2$ be the plane curve corresponding to
$x\in C_d$. $S(0,d)$ is nonsingular at $x$ with tangent
space $H^0(\tilde{X}, \oh_{\proj^2}(d)-N)$ where $N$ is the
divisor of points of $\tilde{X}$ lying over the nodes of $X$.
The additional points correspond to $3d-2$ {\em distinct}
points of $\tilde{X}$. Since $3d-2$ distinct point on $\proj^1$
impose $3d-2$ independent linear
conditions on sections of
$\oh_{\proj^2}(d)-N\eqq \oh_{\proj^1}(3d-2)$, it follows
$C_d$ is nonsingular at $x$.
\end{pf}
Actually, $\barr{M}_{0,0}(\proj^2,d)$ and $S(0,d)$ are
isomorphic on the irreducible, nodal locus. Hence Lemma
(\ref{tt}) is unnecessary.
\begin{lm}
Let $x\in C_d$ be a point corresponding to
a $1$-cuspidal rational plane curve
with all singularities distinct
from $p_1, \ldots, p_{3d-2}$. $C_d$ is cuspidal
at $x$.
\end{lm}
\begin{pf}
The versal deformation space of the cusp $Z_0^2+Z_1^3$ is
2 dimensional:
$$Z_0^2+Z_1^3+ aZ_1 + b.$$
The locus in the versal deformation space corresponding
to equigeneric deformations is determined by
the cuspidal curve $4a^3+27b^2=0$.

Let $X$ be the plane curve corresponding to $x$.
Let $q\in X$ be the cusp. Let
$\tilde{X}$ the normalization of $X$. Let $p\in \tilde{X}$ lie
over $q$.
The nodes of $X$, the
points $p_1, \ldots, p_{3d-2}$, and the $2$ dimensional
subscheme supported on $q$ and pointing in the direction
of the tangent cone of $X$ all together impose
independent conditions on the linear system of degree
$d$ plane curves. First, this independence will be established.

Let $A$ be the subspace of $H^0(\proj^2, \oh_{\proj^2}(d))$
passing through the nodes, points, and the subscheme
of length 2. As before, let $N$ denote the
divisor of $\tilde{X}$ lying above the nodes.
There is a natural left exact sequence obtained by
pulling back sections to $\tilde{X}$:
$$0 \rarr \com \rarr A \rarr
H^0( \tilde{X}, \oh_{\proj^2}(d)-N-p_1-\ldots -p_{3d-2}-3p).$$
By counting conditions,
$$dim(A) \geq {(d+1)(d+2)\over 2} - {(d-1)(d-2)\over 2}+1-
3d+2 -2 = 1$$
with equality if only if the conditions are independent.
Since $$deg_{\tilde{X}}
(\oh_{\proj^2}(d)-N-p_1-\ldots -p_{3d-2}-3p)=
d^2- (d-1)(d-2)+2-3d+2-3 = -1,$$
$dim(A)=1$ and the conditions are independent.

By the independence result above,
the deformations of $X$ parameterized by
the linear system of plane curves through the nodes
and the points $p_1, \ldots, p_{3d-2}$ surjects on the
2 dimensional versal deformation space of the cusp.
The locus of equigeneric deformations of $X$ through the
points $p_1, \ldots, p_{3d-2}$ is \'etale locally
equivalent to the cusp in the versal deformation
space of the cusp.
\end{pf}

\begin{lm}
Let $[\mu]\in \hat{C}_d$ (resp. $x\in C_d$)
be a point corresponding to
an irreducible, nodal, rational curve with a node
at $p_1$ and nodes distinct from $p_2, \ldots, p_{3d-2}$.
$\hat{C}_d$ is nodal at $[\mu]$ (resp. $C_d$ is nodal
at $x$).
\end{lm}
\begin{pf}
If suffices to prove the result for $\hat{C}_d$.
The divisor $D_1 \subset \barr{M}_{0,0}(\proj^2,d)$
corresponding to curves passing
through the point $p_1$ has two nonsingular branches with a normal
crossings intersection at $[\mu]$. Let $r,s\in \proj^1$
lie over $p_1\in \proj^2$.
The two nonsingular branches have the following tangent spaces at $X$:
$$H^0(\proj^1, Norm(-r), \ \ H^0(\proj^1, Norm(-s)).$$
The remaining
$3d-3$ points impose independent conditions on each of these
tangent spaces.
Etale locally at $[\mu]$, $\hat{C}_d$ is the intersection
of the union of linear spaces of dimensions $3d-2$ meeting
along a subspace of dimension $3d-3$ with
general linear space of codimension $3d-3$.
Hence, $\hat{C}_d$ is
nodal at $[\mu]$.
\end{pf}

\begin{lm}
\label{ttt}
Let $x\in C_d$ be a point corresponding to
the union of two irreducible, nodal, rational curves
with degrees $i$ and $d-i$
meeting transversely
with nodes (including component intersection points) distinct from
$p_1, \ldots, p_{3d-2}$. Also assume the components of
degrees $i$, $d-i$ contain $3i-1$, $3(d-i)-1$ points
respectively.
$C_d$ has the singularity type of
the coordinate axes at the origin in $\com^{id-i^2}$.
\end{lm}

\begin{pf}
The nodes (including the intersections of the two components of $X$)
and the points $p_1, \ldots, p_{3d-2}$ necessarily impose
$ (d+1)(d+2)/2$ independent
conditions on $H^0(\proj^2, \oh_{\proj^2}(d))$.
This independence can be established as follows.
Let $\tilde{X}$ be the normalization of $X$ (note $\tilde{X}$ is
the disjoint union of two $\proj^1$'s).
Let $A\subset H^0(\proj^2, \oh_{\proj^2}(d))$ be the linear series
passing through all the nodes of $X$.
There is an exact sequence of vector spaces
$$0 \rarr \com \rarr A \rarr H^0(\tilde{X},
\oh_{\proj^2}(d)-N)\rarr 0.$$
As before, $N$ is the divisor preimage of the nodes of $X$.
Certainly only a $1$ dimensional subspace of $A$
corresponding to the equation of $X$ vanishes on $\tilde{X}$.
Surjectivity of the above sequence follows by a dimension count:
$$dim(A) \geq {(d+1)(d+2)\over 2} - {(d-1)(d-2)\over 2}-1 = 3d-1,$$
$$h^0(\tilde{X}, \oh_{\proj^2}(d)-N)= d^2-(d-1)(d-2)-2+2=3d-2.$$
The points $p_1, \ldots, p_{3d-2}$ are distinct on $\tilde{X}$ and
impose independent conditions on $H^0(\tilde{X}, \oh_{\proj^2}(d)-N)$
by the assumption of their distribution (and the fact $\tilde{X}$ is
a disjoint union of $\proj^1$'s).

At $x\in \barr{S}(0,d)$, the closed Severi variety has
$id-i^2$ nonsingular branches (one for each intersection point).
Let $q\in \proj^2$ be an intersection point of the two components
of $X$. The tangent
space $T(q)$ to the branch of $\barr{S}(0,d)$ corresponding
to $q$ is simply
the linear subspace $T(q)\subset H^0(\proj^2, \oh_{\proj^2}(d))$
of polynomials that vanish at all the nodes of $X$ besides $q$.
Let $V\subset H^0(\proj^2, \oh_{\proj^2}(d))$ be the linear subspace
of polynomials that vanish at the non-intersection nodes of $X$ and
the points $p_1, \ldots, p_{3d-2}$.
$C_d$ at $x$ is \'etale locally equivalent to the intersection
$$V \cap (T(q_1) \cup T(q_2) \cup \cdots \cup T(q_{id-i^2})).$$
Note $V\eqq \com^{id-i^2}$. Since
the nodes of $X$ and the points $p_1, \ldots, p_{3d-2}$ impose
independent conditions on $H^0(\proj^2, \oh_{\proj^2}(d))$, the
Lemma is proven.
\end{pf}

\noindent
The last case to be consider is the analogue of Lemma (\ref{ttt})
for $\hat{C}_d$: when $[\mu]\in \hat{C}_d$ corresponds to
a map with reducible domain and image
satisfying the conditions of (\ref{ttt}). This case can be
handled directly.
However, it is easier to observe that at such $[\mu]$,
$\barr{M}_{0,0}(\proj^2,d)$ is locally isomorphic to the nonsingular
branch in the proof of Lemma (\ref{ttt})
determined by the attaching point of the two components.
The singularity analysis then shows $[\mu]\in \hat{C}_d$ is
a nonsingular point.

For general points $p_1, \ldots, p_{3d-2}$, every point $x\in C_d$,
$[\mu]\in \hat{C}_d$
corresponds to exactly one of the three cases covered by Lemmas (1-5).
Hence the singularities of $C_d$, $\hat{C}_d$ are established.
\label{singa}

\subsection{The Arithmetic Genus}
Consider the moduli space $\barr{M}_{0,0}(\proj^2, d)$ for
$d\geq 3$ (to avoid $[0,0,r,d]=[0,0,2,2]$).
For general points $p_1, \ldots, p_{3d-2}$, the
intersection cycle
$$\hat{C}_d=\HH_1 \cap \HH_2 \cap \cdots \cap \HH_{3d-2}$$
is a curve in $\barr{M}_{0,0}(\proj^2, d)$. $\HH_i$
is the divisor of maps passing through the point $p_i$.
By the analysis in section (\ref{singa}), $\hat{C}_d$ is
nonsingular except for nodes. The nodes occur precisely
at the points $[\mu]\in \hat{C}_d$ corresponding to a
nodal curve with a node at some $p_i$.
Since, for general points,
$\hat{C}_d \subset \barr{M}^*_{0,0}(\proj^2, d)$,
the arithmetic genus
$\hat{g}_d$ of $\hat{C}_d$ can be determined by the
formula for the canonical class
and adjunction ($d\geq 3$):
$$2\hat{g}_d-2 = \big(\KK + (3d-2) \HH\big)\cdot \HH^{3d-2}.$$
A computation of these intersection numbers
in terms of the numbers $N_d$ yields for all $d  \geq 3$:
$$2\hat{g}_d-2=
{(2d-3)(3d+1)\over 2d}N_d +
{1\over 4d}\sum_{i=1}^{d-1} N_i N_{d-i}\Big(3i^2(d-i)^2-4di(d-i)
\Big) {3d-2 \choose 3i-1}.$$
The natural map $\hat{C}_d \rarr C_d$ is
a partial desingularization. The arithmetic  genus
of $C_d$ differs from the arithmetic genus of $\hat{C}_d$
only by the contribution of the singularities of Lemma (3) and (5).
Consider first the cusps in $C_d$ determined by Lemma (3).
The number of these cusps is exactly the number
of $1$-cuspidal, degree $d$, rational curves through $3d-2$ points
in the plane. In [P2] it is shown there are
$${3d-3\over d}N_d + {1\over 2d} \sum _{i=1}^{d-1}
N_iN_{d-i}(3i^2(d-i)^2-2di(d-i)) {3d-2\choose 3i-1}$$
$1$-cuspidal, degree $d$, rational curves through $3d-2$ points.
Each cusp contributes $1$ to the arithmetic genus of $C_d$.
The singularities of Lemma (5) contribute
$${1\over 2} \sum _{i=1}^{d-1}
N_iN_{d-i} (i(d-i)-1)  {3d-2\choose 3i-1}$$
to the arithmetic genus of $C_d$.
The
formula for the arithmetic genus of $C_d$ can
be deduced by adding these contributions to the
formula for $\hat{g}$:

$$2g_d-2=
{6d^2+5d-15\over 2d}N_d +
{1\over 4d}\sum_{i=1}^{d-1} N_i N_{d-i}\Big(15i^2(d-i)^2
-8di(d-i)-4d\Big) {3d-2 \choose 3i-1}.$$

\subsection{The Geometric Genus}
The geometric genus, $g(\tilde{C}_d)$ is simple to compute.
By Bertini's Theorem, the curve $\tilde{C_d}$ determined
in $\barr{M}_{0,3d-2}(\proj^2,d)$ by $3d-2$ general points
is nonsingular and contained in the automorphism-free
locus $\barr{M}^*_{0,3d-2}(\proj^2,d)$.
There is sequence of natural maps exhibiting $\tilde{C}_d$ as
the normalization of both $\hat{C}_d$ and $C_d$:
$$\tilde{C}_d \rarr \hat{C}_d \rarr C_d.$$
The genus of $\tilde{C}_d$
can be determined by the formula for the canonical class
and adjunction:
$$2\tilde{g}_d-2 = \big(\KK+2\sum_{1}^{3d-2} c_1(\LL_p)\big) \cdot
c^2_1(\LL_1) \cdots c^2_1(\LL_{3d-2})$$
$$= \KK \cdot c^2_1(\LL_1) \cdots c^2_1(\LL_{3d-2}).$$
A computation of these intersection numbers in terms of
the numbers $N_d$ yields for all $d\geq 1$:
$$2\tilde{g}_{d}-2 = -{3d^2-3d+4\over 2d^2} N_d
+ {1\over 4d^2} \sum_{i=1}^{d-1} N_i N_{d-i}
(id-i^2) \Big( (9d+4)i(d-i)-6d^2\Big) {3d-2\choose 3i-1}.$$

\subsection{The difference $\hat{g}_d- \tilde{g}_d$.}
Let $d\geq 3$. The natural map $\tilde{C}_d \rarr \hat{C}_d$
is a desingularization. $\hat{C}_d$ has only nodal singularity.
The difference,
$\hat{g}_d - \tilde{g}_d$, equals the number of nodes of
$\hat{C}_d$. Let $M_d$ be
the number of irreducible, nodal, rational
degree $d$ plane curves
with a node at a fixed point and passing through
$3d-3$ general point in $\proj^2$. From the description of the
nodes of $\hat{C}_d$, it follows:
$$\hat{g}_d- \tilde{g}_d= (3d-2) M_d.$$
Values for low degree $d$ are tabulated below:

\begin{tabular}{|l|l|} \hline
$d$ & $N_d$ \\
\hline
1 & 1\\
2 & 1\\
3 & 12\\
4 & 620\\
5 & 87304 \\
6 & 26312976\\
7 & 14616808192 \\
8 & 13525751027392 \\
\hline
\end{tabular}

\begin{tabular}{|l|l|l|l|l|} \hline
$d$& $g_d$&  $\hat{g}_d$&$ \tilde{g}_d$ & $M_d$\\
\hline
1 & 0 & 0 & 0 &  * \\
2 & 0 & 0 & 0 &  * \\
3 & 55 & 10 & 3 & 1 \\
4 & 5447 & 1685 & 725 & 96 \\
5 & 1059729 & 402261 & 166545&  18132 \\
6 & 393308785 & 168879025 & 64776625 &  6506400 \\
7 & 254586817377 & 119342269809 & 42214315809 &  4059366000 \\
8 & 265975021514145 & 133411753757505
& 43616611944513 &  4081597355136 \\
\hline
\end{tabular}

\noindent The formula for $M_d$ (for $d\geq 3$) is:
$$M_d= {d^2-1\over d^2} N_d
- {1\over 4d^2} \sum_{i=1}^{d-1} N_i N_{d-i}
(id-i^2) \Big( {(6d+4)i(d-i)-2d^2\over 3d-2} \Big)
{3d-2\choose 3i-1}.$$
An alternative method of computing $g_3$ is the following.
$\barr{S}(0,3)$ is simply the degree 12 discriminant hypersurface
in the linear system of plane cubics. Therefore, $C_3$ is
a degree 12 plane curve of arithmetic genus $11\cdot10/2=55$.
In fact, $C_3$ has 24 cusp, 28 nodes, and geometric genus 3.

\noindent
Department of Mathematics, University of Chicago

\noindent
rahul@@math.uchicago.edu

\end{document}